\documentclass{amsart}
\usepackage{amsmath}
\usepackage{amsfonts}
\usepackage{amsthm}
\usepackage{amssymb}
\usepackage{enumerate}
\usepackage[all]{xy}
\usepackage{cite}
\usepackage{mathrsfs}
\usepackage{color}
\usepackage{xcolor}
\usepackage{hyperref}
\usepackage{fancyhdr}
\hypersetup{
	colorlinks=true,
	anchorcolor=blue,
	linkcolor=blue,
	filecolor=blue,
	urlcolor=blue,	
	citecolor=blue,
	bookmarks=true,
	bookmarksopen=true,
	pdfborder=000
}
\numberwithin{equation}{section}
\theoremstyle{plain}
\newtheorem{thm}{Theorem}[section]

\newtheorem{proposition}[thm]{Proposition}

\newtheorem{cor}[thm]{Corollary}

\newtheorem{lemma}[thm]{Lemma}
\newtheoremstyle{noparens}%
 {}{}%
 {\itshape}{}%
 {\bfseries}{.}%
 { }%
 {\thmname{#1}\thmnumber{ #2}\mdseries\thmnote{ #3}}
\theoremstyle{noparens}
\newtheorem{lemmaNoParens}[thm]{Lemma}
\newtheorem{thmNoParens}[thm]{Theorem}
\newtheorem{propositionNoParens}[thm]{Proposition}
\theoremstyle{definition}
\newtheorem{defn}[thm]{Definition}
\theoremstyle{remark}
\newtheorem{rmk}[thm]{Remark}

\makeatletter

\newcommand{\Rmnum}[1]{\expandafter\@slowromancap\romannumeral #1@}
\makeatother

\begin{document}

\title{Bi-H\"{o}lder extensions of quasi-isometries on pseudoconvex domains of finite type in $\mathbb{C}^2$}
\author{Jinsong Liu\textsuperscript{1,2} $\&$ Xingsi Pu\textsuperscript{1,2} $\&$ Hongyu Wang\textsuperscript{3}}

\address{$1.$ HLM, Academy of Mathematics and Systems Science,
Chinese Academy of Sciences, Beijing, 100190, China}
\address{$2.$ School of
Mathematical Sciences, University of Chinese Academy of Sciences,
Beijing, 100049, China }
\address{$3.$ School of Science, Beijing University of Posts and Telecommunications, Beijing 100876, China}
\email{liujsong@math.ac.cn,\:puxs@amss.ac.cn,\:wanghyu@bupt.edu.cn}

\begin{abstract}
In this paper, we prove that the identity map for the smoothly bounded pseudoconvex domain of finite type in $\mathbb{C}^2$ extends to a bi-H\"{o}lder map between the Euclidean boundary and Gromov boundary. As an application, we show the bi-H\"{o}lder boundary extensions for quasi-isometries between these domains. Moreover, we get a more accurate index of the Gehring-Hayman type theorem for the bounded $m$-convex domains with Dini-smooth boundary.
\end{abstract}
\maketitle

\section{\noindent{{\bf Introduction}}}
 The boundary extension problems for biholomorphisms and for more generally rough
quasi-isometries play an important role in single or several complex variables.
In the complex plane $\mathbb{C}$ every biholomorphic map between Jordan domains can extend to a homeomorphism between the closures. In $\mathbb{C}^n, n\geq 2$, the boundary extension problem is more complicated.
For example, we have the following classical conclusion from the survey \cite[Theorem 1.1]{1993Proper}:

Suppose that $\Omega_1, \:\Omega_2\subset\mathbb{C}^n(n\geq2)$ are bounded pseudoconvex domains with $C^2$-smooth boundaries, which satisfy
\begin{align}\label{c2}
\frac{1}{C} \delta_{\Omega_1}(z)^{\frac{1}{\alpha}} \leq \delta_{\Omega_2}(f(z)) \leq C \delta_{\Omega_1}(z)^{\alpha}, \quad \forall\; z \in \Omega_1.
\end{align}
And, for any $z \in \Omega_i,\;v \in \mathbb{C}^n$, if there exist contants $\beta_1, \beta_2, C>0$ such that the Kobayashi metric
$$
k_{\Omega_i}(z, v) \geq \frac{C|v|}{\delta_{\Omega_i}(z)^{\beta_i}},
$$
where $\delta_{\Omega_i}(z):=\inf \left\{|w-z|:w \in \partial \Omega_i\right\}, i=1,2$, then every bi-holomorphic map $f: \Omega_1 \rightarrow \Omega_2$ is a bi-H\"{o}lder map.
 Therefore, it can extend to a bi-H\"{o}lder continuous map of $\overline{\Omega}_1$. In particular, the result is valid for strictly pseudoconvex domains or bounded smooth pseudoconvex domains of finite type \cite[Corollary 2]{1992A}. Moreover, on strictly pseudoconvex domains the result is also valid if we substitute quasi-isometries for biholomorphisms\cite[Corollary 6.1]{2000Gromov}.

On the other hand, Balogh and Bonk \cite{2000Gromov} investigated the {\it Gromov hyperbolicity} of strictly pseudoconvex bounded domains in $\mathbb{C}^n$ endowed with Kobayashi metric. And they gave a correspondence between the Euclidean boundary and the Gromov boundary which could be used in boundary extension problems for proper holomorphic mappings or quasi-isometric mappings.

Later Zimmer \cite{Zimmer2016Gromov} and Fiacchi \cite{2022Gromov} proved the Gromov hyperbolicity for bounded convex domains in $\mathbb{C}^n$ or pseudoconvex domain in $\mathbb{C}^2$ with smooth boundary of finite type. They also showed that the Gromov boundary is homeomorphic to the Euclidean boundary which could be used to show the homeomorphic extensions for the quasi-isometric mappings. Furthermore, both the correspondence and the extension results have been improved to be bi-H\"{o}lder for bounded convex domains with smooth boundary of finite type (see \cite[Propsition 12.2]{zimmer22} or \cite[Propsition 1.8]{2022Bi}).

In this paper, we prove that the identity map for the smooth bounded pseudoconvex domain of finite type in $\mathbb{C}^2$ can extend to not only a homeomorphism but also a bi-H\"{o}lder map between the Euclidean boundary and Gromov boundary.

\bigskip
The main theorem is the following:

\begin{thm}\label{id}
Let $\Omega \subset \mathbb{C}^{2}$ be a bounded pseudoconvex domain with smooth boundary of finite D'Angelo type. Then the identity map id : $\Omega \rightarrow \Omega$ can be extended to a bi-H\"{o}lder map
$$
id:(\partial \Omega ,|\cdot|) \rightarrow\left(\partial_G^{d_{K}}\Omega, \rho_G^{d_{K}}\right)
$$
between the boundaries, where $\rho_G^{d_{K}}$ is a visual metric on the Gromov boundary of $\left(\Omega, d_{K}\right)$.
\end{thm}

By using the above theorem, we can generalize the results on bi-H\"{o}lder extensions to the case of pseudoconvex domains with smooth boundary of finite type in $\mathbb{C}^2$.
\begin{thm}\label{bi}
 Let $\Omega_1,\Omega_2 \subset \mathbb{C}^2$ be bounded pseudoconvex domains with smooth boundary of finite D'Angelo type. If $f: \Omega_1 \rightarrow \Omega_2$ is a quasi-isometry with respect to the Kobayashi metrics, then $f$ has a homeomorphic extension $\bar{f}: \overline{\Omega}_1 \rightarrow \overline{\Omega}_2$ such that the induced boundary map $\bar{f}|_{\partial \Omega_1}: \partial \Omega_1 \rightarrow \partial \Omega_2$ is bi-H\"{o}lder continuous with respect to the Euclidean metric.
\end{thm}

Following Balogh and Bonk\cite{2000Gromov}, in order to prove the above theorems, we first give a lower bound for the Finsler distances which satisfy some conditions.

\begin{thm}\label{Fin}
Let $\Omega \subset \mathbb{C}^{n}(n\geq2)$ be a bounded domain and $d_F$ be a complete Finsler distance. Suppose that there exist constants $C_1,C_2,\alpha>0,1>\beta>0$ such that, for all $x,\:y\in\Omega$,
\begin{align}\label{slow}
 d_F(x,y)\geq \alpha\left|\log\left(\frac{\delta(x)}{\delta(y)}\right)\right|-C_1,
\end{align}
and the Finsler metric
\begin{align}\label{fmet}
F(z,X)\geq\frac{C_2\left|X\right|}{\delta(z)^{\beta}}, \:\:\: \forall z\in\Omega,\; 0\neq X\in\mathbb{C}^n.
\end{align}
Then, for each $x,y \in \Omega$,
\begin{align}
d_F(x,y)\geq2\alpha\log\left(\frac{L([x,y])^{\frac{1}{\beta}}+\delta(x)\vee\delta(y)}{\sqrt{\delta(x)\delta(y)}}\right)-C
\end{align}
for some constant $C>0$.

Here $a\vee b:=\max\{a,b\}$ and $L([x,y])$ denotes the Euclidean length of a Finsler geodesic $[x,y]$ connecting $x$ and $y$ in $\Omega$.
\end{thm}

Note that, if the metric is not complete and the conditions are satisfied locally, we can substitute $L(x,y)$ with $|x-y|$.

\bigskip
For the Kobayashi distance on the bounded $m$-convex domains, we have the following result.

\begin{cor}\label{mco}
Let $\Omega \subset \mathbb{C}^{n}$ be a bounded $m$-convex domain. Then there exists a constant $C>0$ such that
$$
d_K(x,y)\geq\log\left(\frac{L([x,y])^m+\delta(x)\vee\delta(y)}{\sqrt{\delta(x)\delta(y)}}\right)-C
$$
for any $x,\: y \in \Omega$.
\end{cor}
\begin{rmk}
For any bounded smooth domain which is locally $\mathbb{C}$-convexifiable near a point of finite type in the boundary, Nikolov and Thomas proved a similar estimate for Kobayashi distance in \cite[Theorem 6]{2016Comparison}. Moreover, Wang has obtained a precise estimate of Kobayashi distance up to a bounded additive term for the bounded convex domains of finite type recently in\cite[Theorem 1.1]{wang22}.
\end{rmk}

Combining the upper bound estimates of Kobayashi distances in Lemma \ref{up} and Corollary \ref{mco}, we can provide another proof of Theorem 1.2 in \cite{2022Bi}.
\begin{thm}\label{miso}
 Let $\Omega_1,\Omega_2 \subset \mathbb{C}^n, n\geq2$, be bounded $m$-convex domains with Dini-smooth boundary, and let $\overline{\Omega}_1,\overline{\Omega}_2$ be their Euclidean compactifications respectively. If $f: \Omega_1 \rightarrow \Omega_2$ is an isometry with respect to the Kobayashi metrics, then $f$ has a homeomorphic extension $\bar{f}: \overline{\Omega}_1 \rightarrow \overline{\Omega}_2$ such that the induced boundary map $\bar{f}|_{\partial \Omega_1}: \partial \Omega_1 \rightarrow \partial \Omega_2$ is bi-H\"{o}lder continuous with respect to the Euclidean metric.
\end{thm}

When these domains are also Gromov hyperbolic, we get the extension results for quasi-isometries. Moreover, we obtain the Gehring-Hayman type Theorem in \cite[Theorem 1.5]{2022Bi} with more accurate index. Refer to \cite[Theorem 14]{NA22} on further details.
\begin{thm}\label{lgh}
Let $\Omega\subset \mathbb{C}^n(n\geq2)$ be a bounded $m$-convex domain with Dini-smooth boundary. Then there exists a constant $C>0$ with
\[L([x,y])\leq C|x-y|^{\frac{1}{m}}\]
for every $x,y\in\Omega$, where $[x,y]$ denotes the Kobayashi geodesic connecting $x$ and $y$.

Moreover, suppose that $(\Omega,d_K)$ is Gromov hyperbolic and $\gamma$ is a Kobayashi $\lambda$-quasi-geodesic connecting $x$ and $y$ with $\lambda\geq1$.
Then there exists a constant $C'>0$ such that
\[L(\gamma)\leq C'|x-y|^{\frac{1}{m}}.\]
\end{thm}
\begin{rmk}
Similarly, the indexes become $1/2$ for strictly pseudoconvex domains. For those bounded smooth pseudoconvex domains of finite type in $\mathbb{C}^2$, we can also get the same results for Kobayashi metrics as they are Gromov hyperbolic and bi-Lipschitz to Catlin metrics.
\end{rmk}

As an application, we obtain a comparision between the local and global Kobayashi distances inspired by Nikolov and Thomas \cite[Theorem 1.1]{2022Qu}.
\begin{thm}\label{add}
Let $\Omega\subset \mathbb{C}^n, n\geq2,$ be a bounded domain. Suppose that $\Omega\cap U$ is $m$-convex and $\partial\Omega\cap U$ is Dini-smooth in a neighborhood $U$ of $\xi\in\partial\Omega$.
Then there exists a neighborhood $V$ of $\xi$ with $V\subset\subset U$ and a constant $C>0$ such that, for any $x,y\in\Omega\cap V$,
\begin{align}\label{loc1}
 d_{K_\Omega}(x,y)\leq d_{K_{\Omega\cap U}}(x,y)\leq d_{K_\Omega}(x,y)+C|x-y|^{\frac{1}{m}},
\end{align}
and
\begin{align}\label{loc2}
1\leq\frac{d_{K_{\Omega\cap U}}(x,y)}{d_{K_\Omega}(x,y)}\leq 1+C\left(\delta_\Omega(x)\wedge\delta_\Omega(y)+|x-y|^{\frac{1}{m}}\right)\text{ for }x\neq y\in\Omega\cap V.
\end{align}
\end{thm}

\begin{rmk}
Nikolov and Thomas \cite[Theorem 1.1]{2022Qu} firstly proved the above localization result for the domain which is smooth and locally $\mathbb{C}$-convexifiable near a point of finite type in the boundary, and Theorem \ref{add} actually relaxes the boundary regularity in the statements of their results.
\end{rmk}

Recall that, for the bounded smooth pseudoconvex domains of finite type in $\mathbb{C}^n$, Cho has proved the condition (\ref{fmet}) in \cite[Theorem 1]{1992A}, but we don't know if they meet the condition (\ref{slow}) when $n\geq3$. Moreover, Herbort \cite[Theorem 2.1]{Herbort2005} has obtained an estimate for the Kobayashi distance by the result in \cite[Theorem 1]{1989Estimates} when $n=2$. In this case, we can get the condition (\ref{slow}), but $\alpha$ may not equal to $\frac{1}{2}$. So in order to prove the corresponding extension result, we have to consider the Catlin metric.

\bigskip
Now suppose that $\Omega =\{z\in\mathbb{C}^{2}:r(z)<0\}$ is a bounded smooth pseudoconvex domain of finite D'Angelo type, and $\xi \in \partial \Omega$ is a point of type $m_\xi$. By using a rotation of the canonical coordinates, we can arrange that the normal direction to $\partial \Omega$ at $\xi$ is given by the $\Im z_1$-axis. Supposing that $\xi=0$, and by using Implicit function theorem, we obtain a local defining function of the form $r\left(z_1, z_2\right)=$ $\Im z_2-F\left(z_1, \Re z_2\right)$, where $F$ is a smooth function and $F(0)=0$. As $\frac{\partial r}{\partial z_2}(\xi) \neq 0$, in a neighborhood of $\xi$ we can define the following vector fields
$$
L_1:=\frac{\partial}{\partial z_1}-\left(\frac{\partial r}{\partial z_2}\right)^{-1} \frac{\partial r}{\partial z_1} \frac{\partial}{\partial z_2}, \text { and } L_2:=\frac{\partial}{\partial z_2}.
$$
Note that $L_1r\equiv0$ and $L_1,L_2$ form a basis of $T_z^{1,0}$ for all $z$ near $\xi$.
For any $j, k>0$, set
$$
\mathcal{L}_{j, k}(z):=\underbrace{L_{1} \ldots L_{1}}_{j-1 \text { times }} \underbrace{\bar{L}_{1} \ldots \bar{L}_{1}}_{k-1 \text { times }} \partial\bar{\partial}r(L_1,\bar{L}_1)(z).
$$

As $\xi \in \partial \Omega$ is a point of type $m_\xi$, from the results in \cite[Theorem 2.4]{1977A}, it follows that there exist $j_0, k_0$ with $j_0+k_0=m_\xi$ which satisfy
$$
\mathcal{L}_{j, k}(\xi) =0 \quad j+k<m_\xi,\text{ and }\mathcal{L}_{j_0, k_0}(\xi) \neq 0.
$$
Denote
$$
C_{l}^{\xi}(z)=\max \left\{\left|\mathcal{L}_{j, k}(z)\right|: j+k=l\right\}.
$$
Let $X=b_1L_{1}+b_2 L_2$ be a holomorphic tangent vector at $z$. Now define the {\it Catlin metric}
$$
M_{\xi}(z, X):=\frac{|b_2|}{|r(z)|}+|b_1| \sum_{l=2}^{m_\xi}\left(\frac{C_{l}^{\xi}(z)}{|r(z)|}\right)^{\frac{1}{l}}.
$$

Owing to the result of Catlin\cite[Theorem 1]{1989Estimates}, the Kobayashi metric is locally bi-Lipschitz to the Catlin metric.
\begin{thm}
Let $\Omega =\{z\in \mathbb{C}^{2}:r(z)<0\}$ be a smoothly bounded pseudoconvex domain. If $\xi \in \partial \Omega$ be a point of finite type $m_\xi$, then there exist a neighborhood $U$ of $\xi$ and a constant $C \geq 1$ such that
$$
\frac{1}{C} M_{\xi}(z, X) \leq K_{\Omega}(z, X) \leq C M_{\xi}(z, X)
$$
for each $z \in \Omega \cap U$ and $X \in \mathbb{C}^{2}$.
\end{thm}

For our convenience, we may make a small change to the form of Catlin metrics. If $X$ is a holomorphic tangent vector at $z$, it has the unique orthogonal decomposition $X=X_{H}+X_{N}$ by $X_{H} \in H_{\pi{(z)}} \partial \Omega$ and $X_{N} \in N_{\pi{(z)}} \partial \Omega$ (Section \ref{boundary}).
As $L_2$ may be not parallel to $X_N$, we set
$$
\widetilde{M}_{\xi}(z, X):=\frac{|X_N|}{\delta(z)}+|X_H| \sum_{l=2}^{m_\xi}\left(\frac{C_{l}^\xi(z)}{\delta(z)}\right)^{\frac{1}{l}}.
$$
Lemma \ref{eq} implies therefore that ${M}_{\xi}(z, X)\asymp\widetilde{M}_{\xi}(z, X)$ in a neighborhood $U$ of $\xi$.

Now choose an open neighborhood $U_i$ of $\xi_i\in\partial\Omega$, $1\leq i\leq s$, which form a finite cover of $\partial\Omega$. There exists a small $\varepsilon>0$ such that the neighborhood
$$
N_{\varepsilon}(\partial \Omega):=\left\{z \in \Omega \| \delta_{\Omega}(z) \mid<\varepsilon\right\}\subset\bigcup_{i=1}^{s}U_i.
$$
Denote $I_z:=\{i:z\in \overline{U}_i\}$, and set
\[\widetilde{M}(z,X):=\max_{i\in I_z}\{
\widetilde{M}_{\xi_i}(z, X)\} \text { for } z\in \Omega\cap N_{\varepsilon}(\partial \Omega).
\]
Since it is upper semi-continuous, we can define a global Finsler metric in $\Omega$ by
$$\widetilde{K}(z,X):=K(z,X)S(z,X)$$
with the positive function $S(z,X)\asymp1$, and $\widetilde{K}(z,X)=\widetilde{M}(z, X)$ for $z\in\Omega\cap N_\varepsilon(\partial\Omega)$,
which implies that $\widetilde{K}(z,X)\asymp K(z,X)$.

Denoting $m=\max\{m_{\xi_i}\}$, for $z\in\Omega\cap N_\varepsilon(\partial\Omega)$ and $X \in \mathbb{C}^{2}$, it now follows that
\begin{align}\label{nk}
\frac{|X_N|}{\delta(z)}+\frac{1}{C}\frac{|X_H|}{\delta(z)^{\frac{1}{m}}}\leq\widetilde{K}(z,X)\leq\frac{|X_N|}{\delta(z)}+C\frac{|X_H|}{\delta(z)^{\frac{1}{2}}}
\end{align}
for some constant $C>0$.
We may call it a Catlin-type metric, and denote by $d_{\widetilde{K}}$ the distance associated to this metric.

\bigskip
Using the technique of Forstneric and Rosay \cite[Proposition 2.5]{1987Localization} and Theorem \ref{Fin}, we get the following estimates.
\begin{thm}\label{dg}
Assume that $\Omega \subset \mathbb{C}^{2}$ is a smoothly bounded pseudoconvex domain of finite type $m$. Then there exists a constant $C \geq 1$ such that
$$
 g_m(x,y)-C\leq d_{\widetilde{K}}(x,y) \leq g_1(x,y)+C
$$
for each $x,y \in \Omega$.
\end{thm}
Here $g_k:\Omega\times\Omega\rightarrow[0,+\infty)$ is defined by
\[g_k(x,y)=2\log\left(\frac{{|x-y|}^{k}+\delta(x)\vee\delta(y)}{\sqrt{\delta(x)\delta(y)}}\right),
\]
and $|x-y|$ is the standard Euclidean distance.

Due to the properties of Gromov hyperbolicity and the estimates above, we deduce that the boundary extensions of quasi-isometries between such domains are bi-H\"{o}lder.

\bigskip
The paper is organized as follows. In Section 2 we will recall some definitions and preliminary results. In Section 3 we give the proof of Theorem \ref{Fin}. In Section 4 we present some applications about boundary extensions of mappings. Finally we will give the proof of Theorem \ref{bi}.

\bigskip

\section{Preliminaries}
\subsection{Notation}
(1) For $z \in \mathbb{C}^n$, let $|\cdot|$ denote the standard Euclidean norm, and let $\left|z_1-z_2\right|$ denote the standard Euclidean distance of $z_1, z_2 \in \mathbb{C}^n$.

(2) Given an open set $\Omega \subsetneq \mathbb{C}^n, x \in \Omega$ and $v \in \mathbb{C}^n$, denote
$$
\delta(x)\text{ or }\delta_{\Omega}(x)=\inf \{|x-\xi|: \xi \in \partial \Omega\}
$$
and
$$\delta_{\Omega}(x,v)=\inf\{|x-\xi|:\xi\in\partial \Omega\cap(x+\mathbb{C}v)\}.$$

(3) Recall that, for real numbers $a, b$, $a \vee b:=\max \{a, b\}$ and $a \wedge b:=\min \{a, b\}$.

(4) For any curve $\gamma$, denote its Euclidean length and Kobayashi length by $L(\gamma)$ and $L_K(\gamma)$ respectively.

(5) For functions $f, g$, write $f \lesssim g$ if there exists $C>0$ such that
$f \leq C g$, and we write
\[f \asymp g\text{ if }f \lesssim g\text{ and }g \lesssim f.\]

\subsection{Boundary projection of curve }\label{boundary}
We first give the relationship between a point and its boundary projection. Note that this result is valid for domains in $\mathbb{R}^{n}(n \geq 2)$ and its proof can be found in \cite[Lemma 2.1]{2000Gromov}.

\begin{lemma}\label{pro}
Suppose that $\Omega=\left\{x \in \mathbb{R}^{n} \mid r(x)<0\right\}, n \geq 2$, is a bounded domain with $C^{2}$-smooth boundary. Then there exists $\delta_{0}>0$ such that

(a). for every point $x \in N_{\delta_{0}}(\partial \Omega)=\left\{z \in \Omega \| \delta_{\Omega}(z) \mid<\delta_{0}\right\}$, there exists a unique point $\pi(x) \in \partial \Omega$ with $|x-\pi(x)|=\delta_{\Omega}(x)$.

(b). the signed distance function $\rho: \mathbb{C}^n \rightarrow \mathbb{R}$
$$
\rho(x)=\left\{\begin{array}{rll}
-\delta(x) & \text { for } & x \in \Omega \\
\delta(x) & \text { for } & x \in \mathbb{C}^n \backslash \Omega
\end{array}\right.
$$
is $C^2$-smooth in $N_{\delta_{0}}(\partial \Omega)$.

(c). for the fibers of the $\operatorname{map} \pi: N_{\delta_{0}}(\partial \Omega) \rightarrow \partial \Omega$, we have
$$
\pi^{-1}(p)=\left(p-\delta_{0} \vec{n}(p), \:\: p+\delta_{0} \vec{n}(p)\right),
$$
where $\vec{n}(p)$ is the outer unit normal vector of $\partial \Omega$ at $p \in \partial \Omega$.

(d). the gradient of the defining function $\rho$ satisfies
$$
\nabla \rho(x)=\vec{n}(\pi(x))
$$
for all $x \in N_{\delta_{0}}(\partial \Omega)$.

(e). the projection map $\pi: N_{\delta_{0}}(\partial \Omega) \rightarrow \partial \Omega$ is $C^{1}$-smooth.
\end{lemma}

Now we start considering domains in complex spaces. It is natural to identify $\left(z_{1}, \ldots, z_{n}\right) \in \mathbb{C}^{n}$ with $ \left(\operatorname{Re} z_{1}, \operatorname{Im} z_{1}, \ldots, \operatorname{Re} z_{n}, \operatorname{Im} z_{n}\right) \in \mathbb{R}^{2n}$.
Recall that, for any $p \in \partial \Omega$, the real tangent space $T_p \partial \Omega$ is given by
$$T_p \partial \Omega=\left\{X \in \mathbb{C}^n: \operatorname{Re}\langle\bar{\partial} r(p), X\rangle=0\right\},$$
and its complex tangent space is
$$H_p \partial \Omega=\left\{X \in \mathbb{C}^n:\langle\bar{\partial} r(p), X\rangle=0\right\},$$
where $$\bar{\partial} r(p)=\left(\frac{\partial r}{\partial \bar{z}_1}(p), \ldots, \frac{\partial r}{\partial \bar{z}_n}(p)\right).$$
Here the standard Hermitian product in $\mathbb{C}^n$ is $\langle X, Y\rangle=\sum_{k=1}^n X_k \overline{Y}_k$. Therefore, for any vector $0\neq X \in \mathbb{C}^{n}$ it has a unique orthogonal decomposition $X=X_{H}+X_{N}$ with $X_{H} \in H_{p} \partial \Omega$ and $X_{N} \in N_{p} \partial \Omega$.

\subsection{The Kobayashi metric}
Given a domain $\Omega\subset\mathbb{C}^n, n\geq2$, a Finsler metric on $\Omega$ is an upper semi-continuous map $F: \Omega \times \mathbb{C}^n \rightarrow [0, \infty)$ with $F(z ; t X)=|t| F(z ; X)$ for any $z \in \Omega, \: t \in \mathbb{C}$ and $X \in \mathbb{C}^n$. The distance function $d_F$ associated with $F$ is defined by
\begin{multline}
\notag d_F(x, y)=\inf \{F\text{-length}(\gamma): \gamma:[0,1] \rightarrow \Omega \text{ is a piecewise }C^1\text{-smooth curve}\\
 \text{with }\gamma(0)=x, \gamma(1)=y\},
\end{multline}
where
$$
F\text{-length}(\gamma)=\int_0^1 F(\gamma(t) ; \dot{\gamma}(t)) d t.
$$

A very important Finsler metric in several complex variables is the Kobayashi metric. For a domain $\Omega \subset \mathbb{C}^n(n \geq 2)$, the (infinitesimal) Kobayashi metric $K_{\Omega}(x ; v), \: x\in\Omega, \: v\in\mathbb{C}^n\backslash\{0\} $, is defined by
$$
K_{\Omega}(x ; v)=\inf \left\{|\xi|: f \in \operatorname{Hol}(\mathbb{D}, \Omega) \text {, with } f(0)=x, d(f)_0(\xi)=v\right\} .
$$
For convenience, we denote by $d_{K_\Omega}$ the Kobayashi distance associated with the Kobayashi metric $K_\Omega$, and sometimes we may omit the subscript $\Omega$. The main property of the Kobayashi distance is that it is contracted by holomorphic maps. That is, if $f: \Omega_1 \rightarrow \Omega_2$ is a holomorphic map, then
$$
\forall z, w \in \Omega_1 \quad d_{K_{\Omega_2}}(f(z), f(w)) \leqslant d_{K_{\Omega_1}}(z, w) .
$$
In particular, the Kobayashi distance is invariant under biholomorphisms, and decreases under inclusions, i.e. if $\Omega_1 \subset \Omega_2 \subset \mathbb{C}^n$ are two bounded domains, then we have $d_{K_{\Omega_2}}(z, w) \leq d_{K_{\Omega_1}}(z, w)$ for all $z, w \in \Omega_1$.

\bigskip
Now we show that the metrics $M_\xi$ and $\widetilde{M}_\xi$ are bi-Lipschitz.
\begin{lemma}\label{eq}
Let $\Omega \subset \mathbb{C}^{2}$ be a smoothly bounded pseudoconvex domain. Supposing that $\xi\in \partial\Omega $ is a boundary point of finite type $m_\xi$, there exist a neighborhood $U$ of $\xi$ and a constant $C \geq 1$ such that, for each $z \in \Omega \cap U$ and $X \in \mathbb{C}^{2}$,
\[\frac{1}{C}{M}_{\xi}(z, X)\leq\widetilde{M}_{\xi}(z, X)\leq C{M}_{\xi}(z, X).
\]
\end{lemma}

\noindent$Proof$. Choose a neighborhood $U$ of $\xi$. Through an affine transformation of the canonical coordinates, without loss of generality we may assume that $\xi=0$ and $\Omega\cap U=\{z\in\mathbb{C}^2:r\left(z_1, z_2\right)=\Im z_2-F\left(z_1, \Re z_2\right)<0\}$, where $F$ is a smooth function such that $F(0)=0$ and $\nabla F(0)=0$. Therefore $\frac{\partial r}{\partial{z_1}}(\xi) =0$, and $\frac{\partial r}{\partial{z_2}}(\xi) \neq 0$.

If $U$ is a chosen sufficiently small neighborhood, then we can define a function $s:=(\frac{\partial r}{\partial{z_2}})^{-1}\frac{\partial r}{\partial{z_1}}$ and have the inequality $0\leq\left|s\right|\leq\frac{1}{2}$ on $U$.
Noting that the vector fields
$$
L_1=\frac{\partial}{\partial z_1}-s\frac{\partial }{\partial z_2}, \text { and } L_2=\frac{\partial}{\partial z_2}
$$
form a basis of $T_z^{1,0}$ for all $z$ near $\xi$, then any $X\in T_z^{1,0}\cong\mathbb{C}^2$ can be written as $X=b_1L_1+b_2L_2=X_H+X_N$, where $X_H\in H_{\pi(z)}(\partial\Omega)$ and $X_N\in N_{\pi(z)}(\partial\Omega)$.

$L_1r\equiv0$ implies that $L_1$ is parallel to $X_H$. Noting that $L_2$ may be not orthogonal to $X_H$, we denote
$$L_1^{\perp}:=s\frac{\partial }{\partial z_1}+\frac{\partial}{\partial z_2},$$
which is parallel to $X_N$.
Therefore,
\begin{align}
\notag b_1L_1+b_2L_2&=b_1L_1+b_2\left(L_1^{\perp}-sL_1\right)\frac{1}{1+s^2}\\
\notag &=\left(b_1-\frac{b_2 s}{1+s^2}\right)L_1+\frac{b_2}{1+s^2}L_1^{\perp},
\end{align}
which implies
\[X_H=\left(b_1-\frac{b_2 s}{1+s^2}\right)L_1,\;\;X_N=\frac{b_2}{1+s^2}L_1^{\perp}
\]
and
\[\left|X_H\right|=\left|b_1-\frac{b_2 s}{1+s^2}
\right|\sqrt{1+\left|s\right|^2},\;\;
\left|X_N\right|=\left|b_2\right|\frac{\sqrt{1+\left|s\right|^2}}{\left|1+s^2\right|}.
\]

For any $z\in\Omega\cap U$, denoting $\omega=\pi(z)\in\partial\Omega\cap U,$ then $\Im \omega_2=F\left(\omega_1, \Re \omega_2\right)$. By the choice of the defining function, we have $$|r(z)|\asymp \delta(z).$$ Then, in view of the above equations, it follows that
\begin{align}
\notag \widetilde{M}_{\xi}(z, X)&=\frac{|X_N|}{\delta(z)}+|X_H| \sum_{l=2}^{m_\xi}\left(\frac{C_{l}^{\xi}(z)}{\delta(z)}\right)^{\frac{1}{l}}
\asymp\frac{|X_N|}{|r(z)|}+|X_H| \sum_{l=2}^{m_\xi}\left(\frac{C_{l}^{\xi}(z)}{|r(z)|}\right)^{\frac{1}{l}}\\
\notag &\asymp\frac{|b_2|}{|r(z)|}+\left|b_1+b_1s^2-b_2s
\right| \sum_{l=2}^{m_\xi}\left(\frac{C_{l}^{\xi}(z)}{|r(z)|}\right)^{\frac{1}{l}}\\
\notag &\asymp\frac{|b_2|}{|r(z)|}+|b_1| \sum_{l=2}^{m_\xi}\left(\frac{C_{l}^{\xi}(z)}{|r(z)|}\right)^{\frac{1}{l}}
\pm|b_2|\left|s\right|
\sum_{l=2}^{m_\xi}\left(\frac{C_{l}^{\xi}(z)}{|r(z)|}\right)^{\frac{1}{l}}\\
\notag &\asymp M_{\xi}(z, X).
\end{align}
The last inequality holds because
\[\frac{1}{|r(z)|}
\pm\left|s\right|
\sum_{l=2}^{m_\xi}\left(\frac{C_{l}^{\xi}(z)}{|r(z)|}\right)^{\frac{1}{l}}\asymp\frac{1}{|r(z)|}.
\]
$\hfill\qed$

Now we give some estimates of the Catlin-type distances on pseudoconvex domains of finite type in $\mathbb{C}^2$.
\begin{lemma}\label{low}
Let $\Omega \subset \mathbb{C}^{2}$ be a smoothly bounded pseudoconvex domain of finite type $m$. Then there exist constants $\varepsilon>0$ and $C \geq 1$ such that, for any $x,\: y\in \Omega\cap N_\varepsilon(\partial\Omega)$,
$$
d_{\widetilde{K}}(x,y) \geq\left|\log \left(\frac{\delta(y)}{\delta(x)}\right)\right|.
$$
Moreover, when $\pi(x)=\pi(y)$, we have
$$
d_{\widetilde{K}}(x,y) =\left|\log \left(\frac{\delta(y)}{\delta(x)}\right)\right|.
$$
\end{lemma}

\noindent$Proof$. Suppose $\gamma:[0,1] \rightarrow \Omega$ is a piecewise $C^1$-smooth curve with endpoints $x=\gamma(0)$ and $y=\gamma(1)$. If $\gamma\subset N_\varepsilon(\partial\Omega)$, by the conclusion (d) in Lemma \ref{pro}, it follows that $2\bar{\partial} \delta\left(\gamma(t)\right)=\nabla\rho\left(\gamma(t)\right)=\vec{n}(\pi(\gamma(t)))$. Then for those $t \in[0,1]$ where $\dot{\gamma}\left(t\right)$ exists, we have
\begin{align}
\notag\left|\frac{d}{d t} \delta(\gamma(t))\right|&=\left|\sum_{k=1}^n\dot{\gamma}(t)\frac{\partial}{\partial z_k}\delta(\gamma(t))+\sum_{k=1}^n\overline{\dot{\gamma}(t)}\frac{\partial}{\partial \bar{z}_k}\delta(\gamma(t))\right|\\
\notag&=2\left|\operatorname{Re}\left\langle\bar{\partial} \delta\left(\gamma(t)\right),\dot{\gamma}(t)\right\rangle\right|
=\left|\operatorname{Re}\left\langle\vec{n}(\pi(\gamma(t))),\dot{\gamma}(t)\right\rangle\right|\\
\notag&\leq \left|\dot{\gamma}_N(t)\right|.
\end{align}
Then (\ref{nk}) implies that
\[\widetilde{K}\text {-length }(\gamma) \geq \int_0^1\frac{|\dot{\gamma}_N(t)|}{\delta(\gamma(t))}d t
\geq \left|\int_0^1 \frac{d(\delta(\gamma(t)))}{\delta(\gamma(t))}\right|
\geq\left|\log \left(\frac{\delta(y)}{\delta(x)}\right)\right|.
\]
If $\gamma$ is not in $N_\varepsilon(\partial\Omega)$, we may assume $\delta(x)\leq\delta(y)$, and set
\[t_0:=\min\{t\in[0,1]:\gamma(t)\notin N_\varepsilon(\partial\Omega)\}.\]
Then $\gamma_0:=\gamma|_{[0,t_0]}\subset N_\varepsilon(\partial\Omega)$, which implies that
\[\widetilde{K}\text {-length }(\gamma) \geq \widetilde{K}\text {-length }(\gamma_0) \geq \log \left(\frac{\varepsilon}{\delta(x)}\right)\geq\log \left(\frac{\delta(y)}{\delta(x)}\right).\]
As $\gamma$ is arbitrary, we deduce that
$$
d_{\widetilde{K}}(x,y) \geq\left|\log \left(\frac{\delta(y)}{\delta(x)}\right)\right|.
$$

Suppose that $\pi(x)=\pi(y)=p\in \partial \Omega$. Since $\gamma:[0,1] \rightarrow \Omega\cap N_\varepsilon(\partial\Omega), t \mapsto x+t(y-x)$ is a straight line segment contained in the fiber $\pi^{-1}(p)$, it follows from Lemma \ref{pro} that
$$
\dot{\gamma}(t) \equiv \pm|x-y| \vec{n}(p), \:\:\: \vec{n}(\pi(\gamma(t)))=\vec{n}(p), \:\:\: \forall t \in[0,1].
$$
So $\frac{d}{d t} \delta(\gamma(t))=\pm|x-y|=\pm\left|\dot{\gamma}_N(t)\right|$, and $\dot{\gamma}_H(t) \equiv 0$.
A similar computation yields
$$
\widetilde{K}\text{-length }(\gamma) =\left|\log \left(\frac{\delta(y)}{\delta(x)}\right)\right|,
$$
which implies
$$
d_{\widetilde{K}}(x,y)=\left|\log \left(\frac{\delta(y)}{\delta(x)}\right)\right|.
$$

In this case, we know that the straight line segment is a geodesic.
$\hfill\qed$

\subsection{The $m$-convex domains}\label{convex}

We first recall the definition of $m$-convex. It could be found in \cite{1993Complex}.
\begin{defn}
A bounded convex domain $\Omega \subset \mathbb{C}^n, n \geq 2$, is called $m$-convex if for some $C>0$ and $m>1$
$$
\delta_{\Omega}(z ; v) \leq C \delta_{\Omega}(z)^{\frac{1}{m}}, \quad z \in \Omega, v \in \mathbb{C}^n.
$$
\end{defn}

For convex domains, Kobayashi metrics have the following estimate.
\begin{lemmaNoParens}[{\cite[Theorem 5]{2017boundary}}]\label{con}
Let $\Omega \subset \mathbb{C}^n$ be a bounded convex domain. Then the Kobayashi metric satisfies
$$
\frac{|v|}{2 \delta_{\Omega}(p ; v)} \leq K_\Omega(p ; v) \leq \frac{|v|}{\delta_{\Omega}(p ; v)}, \:\: \text{\;\; for any }p\in \Omega,\;v \in \mathbb{C}^n.
$$
\end{lemmaNoParens}
\begin{lemmaNoParens}[{\cite[Proposition 2.4]{1993Complex}}]\label{lowconvex}
Let $\Omega \subset \mathrm{C}^n$ be a bounded convex domain. Then the Kobayashi distance
$$
d_K(x, y) \geq \frac{1}{2}\left|\log \frac{\delta_{\Omega}(x)}{\delta_{\Omega}(y)}\right|, \:\:\: \text{\;\; for any } x, y\in \Omega.
$$
\end{lemmaNoParens}

Recall that a $C^1$-smooth boundary point $p$ of a domain $\Omega$ in $\mathbb{C}^n$ is said to be Dini-smooth, if the outer unit normal vector $\vec{n}$ to $\partial \Omega$ near $p$ is a Dini-continuous function. This means that there exists a neighborhood $U$ of $p$ with
$$
\int_0^1 \frac{\omega(t)}{t} d t<+\infty,
$$
where
$$
\omega(t)=\omega(\vec{n}, \partial \Omega \cap U, t):=\sup \left\{\left|\vec{n}(x)-\vec{n}(y)\right|:|x-y|<t, \:\: x, y \in \partial \Omega \cap U\right\}
$$
is the respective modulus of continuity. Note that Dini-smooth is a weaker condition than $C^{1, \epsilon}$-smooth. Here a Dini-smooth domain means that each boundary point of $\Omega$ is a Dini-smooth point.

Then we have the following upper bound of Kobayashi distance.
\begin{lemmaNoParens}[{\cite[Corollary 8]{2015Estimates}}]\label{up}
Let $\Omega$ be a Dini-smooth bounded domain in $\mathbb{C}^n$ and $x, y \in \Omega$.
Then there exists a constant $A>1+\sqrt{2} / 2$ with
$$
d_K(x, y) \leq \log \left(1+\frac{A|x-y|}{\sqrt{\delta_{\Omega}(x) \delta_{\Omega}(y)}}\right).
$$
\end{lemmaNoParens}

\subsection{Gromov hyperbolic}\label{gh}

In this section we will give some definitions and results about Gromov hyperbolicity. Refer to \cite{Metric99} for further details.

Let $(X, d)$ be a metric space. The $Gromov\text{ } product$ of two points $x, y \in X$ with respect to a base point $\omega \in X$ is defined by
$$(x|y)_\omega:=\frac{1}{2}\left(d(x, \omega)+d(y, \omega)-d(x, y)\right).$$

Recall that a metric space $(X, d)$ is a $geodesic$ space if any two distinct points $x, y\in X$ can be joined by a geodesic segment, i.e., the image of an isometry $\gamma:[0, d(x, y)] \rightarrow X$ with $\gamma(0)=x$ and $\gamma(d(x, y))=y$. In particular, if $\gamma$ is a $(\lambda,k)$-quasi-isometry, then we call it a $(\lambda,k)$-quasi-geodesic. And we call it a $\lambda$-quasi-geodesic when $k=0$.

A metric space $(X, d)$ is called $proper$ if every closed ball in $(X, d)$ is compact.
A proper geodesic metric space $X$ is called $Gromov \text{ }hyperbolic$ if there is a constant $\delta\geq0$ such that, for any $x, y, z, \omega \in X$,
$$
(x|y)_\omega \geq \min \left\{(x | z)_\omega,(z | y)_\omega\right\}-\delta.
$$
Then we recall the Gromov boundary.
\begin{defn}
(1) A sequence $\left\{x_i\right\}$ in $X$ is called a $Gromov\text{ } sequence$ if $\left(x_i |x_j\right)_\omega \rightarrow \infty$ as $i$, $j \rightarrow \infty$.

(2) Two such sequences $\left\{x_i\right\}$ and $\left\{y_j\right\}$ are said to be $equivalent$ if $\left(x_i | y_i\right)_\omega \rightarrow \infty$ as $i \rightarrow \infty$.

(3) The $Gromov \text{ }boundary$ $\partial_G X$ of $X$ is defined to be the set of all equivalence classes of Gromov sequences, and $\overline{X}^G=X \cup \partial_G X$ is called the $Gromov \text{ }closure$ of $X$.

(4) For $a \in X$ and $b \in \partial_G X$, the Gromov product $(a | b)_\omega$ of $a$ and $b$ is defined by
$$
(a| b)_\omega=\inf \left\{\liminf _{i \rightarrow \infty}\left(a | b_i\right)_\omega:\left\{b_i\right\} \in b\right\} .
$$

(5) For $a, b \in \partial_G X$, the Gromov product $(a|b)_\omega$ of $a$ and $b$ is defined by
$$
(a|b)_\omega=\inf \left\{\liminf _{i \rightarrow \infty}\left(a_i | b_i\right)_\omega:\left\{a_i\right\} \in a \text { and }\left\{b_i\right\} \in b\right\} .
$$
\end{defn}

For a Gromov hyperbolic space $X$, one can define a class of visual metrics on $\partial_G X$ via the extended Gromov products. For any metric $\rho_G$ in this class, there exist a parameter $\epsilon>0$ and a base point $\omega\in X$ with
$$
\rho_G(a, b) \asymp \exp \left(-\epsilon(a|b)_w\right), \quad \text { for } a, b \in \partial_G X .
$$

Now we recall the definition of rough quasi-isometric maps as follows.

\begin{defn}
Let $f: X \rightarrow Y$ be a map between metric spaces $X$ and $Y$.

(1) If for all $x, y \in X$, there exist constants $\lambda \geq 1,k \geq 0$ with
$$
\lambda^{-1} d_X(x, y)-k \leq d_Y(f(x), f(y)) \leq \lambda d_X(x, y)+k,
$$
then $f$ is called a $(\lambda, k)$-quasi-isometry. In particular, $f$ is called an isometry when $\lambda=1$ and $k=0$.

(2)
For a bijection $f$ and all $x, y \in X$, if there exist constants $\lambda \geq 1,1>\alpha>0$ with
$$
\lambda^{-1}d_X(x, y)^{\frac{1}{\alpha}}\leq d_Y(f(x), f(y))\leq\lambda d_X(x, y)^{\alpha}
$$
then $f$ is $\alpha$-bi-H\"{o}lder.
\end{defn}

Then from the Proposition 5.5 in \cite{Bonk2000}, we have the following extension result:
\begin{proposition}\label{extend}
Let $X$ and $Y$ be Gromov hyperbolic metric spaces. If $f: X \rightarrow Y$ is a quasi-isometry, then $f$ induces a bi-H\"{o}lder map $\tilde{f}: \partial_G X \rightarrow \partial_G Y$.
\end{proposition}

By the stability of quasi-geodesics, we deduce that the Gromov hyperbolicity is quasi-isometric invariant.

\begin{thmNoParens}[{\cite[Prat \Rmnum{3}: Theorem 1.7]{Metric99}}]\label{stab}
Suppose that $(X,d)$ is a $\delta$-Gromov hyperbolic geodesic space with $\delta>0$. And suppose that $\gamma$ is a $(\lambda,\varepsilon)$-quasi-geodesic in $X$ and $[p, q]$ is a corresponding geodesic segment joining the endpoints of $\gamma$. Then there exists a constant $R = R(\delta,\lambda,\varepsilon)$ which satisfies that the Hausdorff distance between $[p, q]$ and the image of $\gamma$ is less than $R$.
\end{thmNoParens}

\begin{propositionNoParens}[{\cite[Prat \Rmnum{3}: Theorem 1.9]{Metric99}}]\label{invariant}
Let $X$ and $Y$ be geodesic metric spaces and let $f: X\rightarrow Y$ be a quasi-isometric embedding. If $Y$ is Gromov hyperbolic, then $X$ is also Gromov hyperbolic.
\end{propositionNoParens}

Given a smooth bounded pseudoconvex domain $\Omega \subset \mathbb{C}^{2}$ of finite type $m$, if $\varepsilon$ is sufficiently small and $\omega\in \Omega\cap N_\varepsilon(\partial\Omega)$ is fixed,
Lemma \ref{low} implies that the distance $d_{\widetilde{K}}(\omega,x)\rightarrow \infty$ as $x\in \Omega \cap N_\varepsilon(\partial\Omega)$ tends to $\partial\Omega$. It implies that the metrics $d_K$ and $d_{\widetilde{K}}$ are both complete.

The Hopf-Rinow Theorem (see {\cite[Prat \Rmnum{1}: Proposition 3.7]{Metric99}}) tells us that any two points in $\Omega$ can be connected by a Kobayashi or Catlin-type geodesic.
In particular, the bounded convex domains endowed the Kobayashi metric are also geodesic spaces \cite{1980Convex}.
$\Omega$ is Gromov hyperbolic with Kobayashi metric by Fiacchi \cite[Theorem 1.1]{2022Gromov}.

In view of Proposition \ref{extend} and \ref{invariant}, we have the following key result.
\begin{proposition}\label{visual}
Let $\Omega\subset \mathbb{C}^2$ be a bounded pseudoconvex domain with smooth boundary of finite type. Then the identity map id : $\Omega \rightarrow \Omega$ can extend to a bi-H\"{o}lder map
$$
id:\left(\partial_G^{d_{\widetilde{K}}}\Omega, \rho_G^{d_{\widetilde{K}}}\right) \rightarrow\left(\partial_G^{d_{K}}\Omega, \rho_G^{d_{K}}\right)
$$
between the boundaries, where $\rho_G^{d_{\widetilde{K}}}$ and $\rho_G^{d_{K}}$ are the visual metrics on the Gromov boundary of $\left(\Omega, d_{\widetilde{K}}\right)$ and $\left(\Omega, d_{K}\right)$ respectively.
\end{proposition}

\bigskip
\section{Estimate of the Finsler metric}

In this section we first prove Theorem \ref{Fin}. For convenience, we restate it as follows.
\begin{thm}
Let $\Omega \subset \mathbb{C}^{n}(n\geq2)$ be a bounded domain and let $d_F$ be a complete Finsler distance. Suppose that there exist constants $C_1,C_2,\alpha>0,1>\beta>0$ with
\begin{align}\label{slow2}
 d_F(x,y)\geq \alpha\left|\log\left(\frac{\delta(x)}{\delta(y)}\right)\right|-C_1, \:\: \text{ for any }x,y\in\Omega
\end{align}
and the Finsler metric
\begin{align}\label{fmet2}
F(z,X)\geq\frac{C_2\left|X\right|}{\delta(z)^{\beta}}, \:\: \text{ for any } z\in\Omega,\; 0\neq X\in\mathbb{C}^n.
\end{align}
Then, for each $x,y \in \Omega$, it follows
\begin{align}
d_F(x,y)\geq2\alpha\log\left(\frac{L([x,y])^{\frac{1}{\beta}}+\delta(x)\vee\delta(y)}{\sqrt{\delta(x)\delta(y)}}\right)-C,
\end{align}
for some constant $C>0$. Here $L([x,y])$ denotes the Euclidean length of a Finsler geodesic $[x,y]$ connecting $x$ and $y$ in $\Omega$.
\end{thm}

\noindent$Proof$.
For any $x,y \in \Omega$, the hypothesis (\ref{slow2}) implies that
$$
d_F(x,y)\geq 2\alpha\log\left(\frac{\delta(x)\vee\delta(y)}{\sqrt{\delta(x)\delta(y)}}\right)-C.
$$
Therefore, we only need to check that
$$
d_F(x,y)\geq 2\alpha\log \left(\frac{L([x,y])^{\frac{1}{\beta}}}{\sqrt{\delta(x)\delta(y)}}\right)-C.
$$

Now for any Finsler geodesic $\gamma:[0,1] \rightarrow \Omega$ with $\gamma(0)=x$ and $\gamma(1)=y$, we know that $\gamma$ is absolutely continuous. Let $L([x,y])$ denote the Euclidean length of $\gamma=[x,y]$. Define $H:=\max _{z \in \gamma} \delta(z)$. There exists $t_0 \in[0,1]$ with $H=\delta\left(\gamma(t_0)\right)$.
Considering the subcurves $\gamma_1=\gamma |_{[0, t_0]}$ and $\gamma_2=\gamma |_{[t_0, 1]}$, there are two possibilities:

If $H \geq L([x,y])^{\frac{1}{\beta}}$, it follows from (\ref{slow2}) that
$$
F \text {-length }\left(\gamma_1\right) \geq \alpha\log \left(\frac{H}{\delta(x)}\right),
$$
and
$$
F\text{-}\operatorname{length}\left(\gamma_2\right) \geq \alpha\log \left(\frac{H}{\delta(y)}\right).
$$
Thus
$$
F\text {-length }(\gamma) \geq 2\alpha\log \left(\frac{H}{\sqrt{\delta(x)\delta(y)}}\right)\geq 2\alpha\log \left(\frac{L([x,y])^{\frac{1}{\beta}}}{\sqrt{\delta(x)\delta(y)}}\right),
$$
which completes the proof.

The other possibility is $H<L([x,y])^{\frac{1}{\beta}}$. Since $\delta(x) \leq H$, there exists $k \in \mathbb{N}_+$ with
$$2^{-\frac{k}{\beta}} H<\delta(x) \leq 2^{-\frac{k-1}{\beta}} H.$$
Consider the curve $\gamma_1$ and define $0=s_0 \leq s_1<\cdots<s_k \leq t_0$ as follows,
$$s_j=\min \left\{s \in\left[0, t_0\right]: \delta(\gamma(s))=\frac{H}{2^{\frac{k-j}{\beta}}}\right\}, \:\:\: j=1, \ldots, k.$$
By denoting $x_j=\gamma\left(s_j\right), \: j=0, \ldots, k$, we have
$$1 \leq \frac{\delta\left(x_j\right)}{\delta\left(x_{j-1}\right)} \leq 2^{\frac{1}{\beta}}.$$

Then we shall consider the following two alternatives:

(a). In the first case we assume that there exists an index $l \in\{1, \ldots, k\}$ with
$$
L([x_{l-1},x_l])>\frac{1}{8} 2^{-(k-l)} L([x,y]).
$$
For convenience, here $[x_{l-1},x_l]$ means $\gamma|_{[s_{l-1}, s_l]}$.
Then, for $t \in\left[s_{l-1}, s_l\right]$, we have
$$
\delta(\gamma(t)) \leq 2^{-\frac{k-l}{\beta}} H,
$$
which implies that
\begin{align}
\notag F\text{-}\operatorname{length}\left(\gamma |_{\left[s_{l-1}, s_l\right]}\right)
& \gtrsim \int_{s_{l-1}}^{s_l}\frac{|\dot{\gamma}(t)|}{{\delta(\gamma(t))}^{\beta}} dt
\geq\frac{2^{k-l}}{H^{\beta}} \int_{s_{l-1}}^{s_l}|\dot{\gamma}(t)| dt\\
&=\frac{2^{k-l}L([x_{l-1},x_l])}{H^{\beta}}
\notag\gtrsim\frac{L([x,y])}{H^{\beta}}.
\end{align}

Let $t_1:=s_k \leq t_0$. From (\ref{slow2}) it follows that
$$
\begin{aligned}
F \text {-length }\left(\gamma|_{\left[0, t_1\right]}\right) &=F\text {-length }\left(\gamma|_{\left[0, s_{l-1}\right]}\right)+F \text {-length }\left(\gamma|_{\left[s_{l-1}, s_l\right]}\right) \\
&+F\text {-length }\left(\gamma|_{\left[s_l, s_k\right]}\right) \\
& \geq \alpha\log \left(\frac{\delta\left(x_{l-1}\right)}{\delta\left(x_0\right)}\right)+C \frac{L([x,y])}{H^{\beta}}+\alpha\log \left(\frac{\delta\left(x_k\right)}{\delta\left(x_l\right)}\right)\\
&\geq\alpha\log \left(\frac{\delta\left(x_k\right)}{\delta\left(x_0\right)}\right)+C \frac{L([x,y])}{H^{\beta}}-\frac{\alpha}{\beta}\log2\\
&=\alpha\log \left(\frac{H}{\delta(x)}\right)+C\frac{ L([x,y])}{H^{\beta}}-C.
\end{aligned}
$$

(b). The second alternative is
$$
L([x_{j-1},x_j]) \leq \frac{1}{8} 2^{-(k-j)} L([x,y]), \:\:\: j=1, \ldots, k.
$$
Then
$$
L([x,\gamma\left(t_1\right)])=\sum_{j=1}^k L([x_{j-1},x_j]) \leq \frac{1}{4} L([x,y]).
$$
follows. On the other hand, by using (\ref{slow2}) once more, we have thus proved that
$$
F\text {-length }\left(\gamma |_{\left[0, t_1\right]}\right) \geq \alpha\log \left(\frac{H}{\delta(x)}\right).
$$
By summarizing this discussion, we obtain the following two possibilities
\begin{equation}\label{e3}
F \text {-length }\left(\gamma|_{\left[0, t_1\right]}\right) \geq \alpha\log \left(\frac{H}{\delta(x)}\right)+C\frac{L([x,y])}{H^{\beta}}-C,
\end{equation}
or
\begin{equation}\label{e4}
F \text {-length }\left(\gamma |_{\left[0, t_1\right]}\right) \geq \alpha\log \left(\frac{H}{\delta(x)}\right), \text { and } L([x,\gamma\left(t_1\right)]) \leq \frac{1}{4} L([x,y]),
\end{equation}
where $t_1 \in\left[0, t_0\right]$.

By applying similar considerations to the curve $\gamma_2$ instead of $\gamma_1$, we can find $t_2 \in\left[t_0, 1\right]$ such that one of the following alternatives is valid
\begin{equation}\label{e5}
F\text {-length }\left(\gamma|_{\left[t_2,1\right]}\right) \geq \alpha\log \left(\frac{H}{\delta(y)}\right)+C\frac{L([x,y])}{H^{\beta}}-C,
\end{equation}
or
\begin{equation}\label{e6}
F \text {-length }\left(\gamma |_{\left[t_2,1\right]}\right) \geq \alpha\log \left(\frac{H}{\delta(y)}\right), \text { and } L([y,\gamma\left(t_2\right)]) \leq \frac{1}{4} L([x,y]).
\end{equation}
Combining (\ref{e4}) with (\ref{e6}), the inequality
$$
L([\gamma\left(t_1\right),\gamma\left(t_2\right)])=L([x,y])-L([x,\gamma\left(t_1\right)])-L([y,\gamma\left(t_2\right)]) \geq \frac{1}{2} L([x,y]),
$$
follows. It is analogous to the estimation in the first alternative. Therefore, we now obtain
$$
F\text {-length }\left(\gamma |_{\left[t_1, t_2\right]}\right) \geq C\frac{L([x,y])}{H^{\beta}}.
$$
Consequently,
$$
\begin{aligned}
F \text {-length }(\gamma) &=F \text {-length }\left(\gamma |_{\left[0, t_1\right]}\right)+F \text {-length }\left(\gamma|_{\left[t_1, t_2\right]}\right)\\
&+F\text {-length }\left(\gamma |_{\left[t_2, 1\right]}\right)
\geq 2\alpha\log \left(\frac{H}{\sqrt{\delta(x) \delta(y)}}\right)+C\frac{L([x,y])}{H^{\beta}}-C.
\end{aligned}
$$
This inequality is also true if (\ref{e3}) and (\ref{e5}) or (\ref{e3}) and (\ref{e6}), or (\ref{e4}) and (\ref{e5}) hold simultaneously. More explicitly,
$$
\begin{aligned}
F \text {-length }(\gamma) &\geq F \text {-length }\left(\gamma |_{\left[0, t_1\right]}\right)+F\text {-length }\left(\gamma |_{\left[t_2, 1\right]}\right)\\
&\geq 2\alpha\log \left(\frac{H}{\sqrt{\delta(x) \delta(y)}}\right)+C\frac{L([x,y])}{H^{\beta}}-C,
\end{aligned}
$$
which implies that the above estimate is true in any case.

Let
$$
f(t):=2\alpha\log \frac{t}{\sqrt{\delta(x) \delta(y)}}+C\frac{L([x,y])}{t^{\beta}}.
$$
Through a simple calculation, we know the function $f$ gets its minimum value when
$$
t=\left(\frac{\beta C}{2\alpha}\right)^{\frac{1}{\beta}}L([x,y])^{\frac{1}{\beta}}
$$
which gives the lower bound
$$
F \text {-length }(\gamma) \geq 2\alpha\log \left(\frac{L([x,y])^{\frac{1}{\beta}}}{\sqrt{\delta(x)\delta(y)}}\right)-C.
$$
As $\gamma$ is a Finsler geodesic, then
$$
d_F(x,y) \geq 2\alpha\log \left(\frac{L([x,y])^{\frac{1}{\beta}}}{\sqrt{\delta(x)\delta(y)}}\right)-C,
$$
which completes the proof.
$\hfill\qed$

Then by using Lemma \ref{con} and Lemma \ref{lowconvex} in Section \ref{convex}, we have the following result.
\begin{cor}\label{mco2}
Let $\Omega \subset \mathbb{C}^{n}$ be a bounded $m$-convex domain. There exists a constant $C>0$ such that the Kobayashi distance
$$
d_K(x,y)\geq\log\left(\frac{L([x,y])^m+\delta(x)\vee\delta(y)}{\sqrt{\delta(x)\delta(y)}}\right)-C
$$
for each $x,y \in \Omega$.
\end{cor}

Note that Lemma \ref{up} implies the Gehring-Hayman theorem for bounded $m$-convex domains.
\begin{thm}
Suppose that $\Omega\subset \mathbb{C}^n(n\geq2)$ is a bounded $m$-convex domains with Dini-smooth boundary. Then there exists a constant $C>0$ such that, for every $x,y\in\Omega$,
\[L([x,y])\leq C|x-y|^{\frac{1}{m}}\]
where $[x,y]$ is a Kobayashi geodesic connecting $x$ and $y$.

Moreover, if $(\Omega,\: d_K)$ is Gromov hyperbolic and $\gamma$ is a Kobayashi $\lambda$-quasi-geodesic connecting $x$ and $y$ with $\lambda\geq1$,
then there exists a constant $C'>0$ which satisfies
\[L(\gamma)\leq C'|x-y|^{\frac{1}{m}}.\]
\end{thm}
\noindent$Proof$. For $x,y\in\Omega$, from Lemma \ref{up} and Corollary \ref{mco2}, it follows that
\[L([x,y])^m\lesssim |x-y|+\sqrt{\delta(x)\delta(y)}.\]
As a consequence, if $\sqrt{\delta(x)\delta(y)}\leq|x-y|$, we have the desired estimation $L([x,y])\lesssim |x-y|^{\frac{1}{m}}$.

When $\sqrt{\delta(x)\delta(y)}>|x-y|$, denoting the curve by $\alpha=[x,y]$ and putting $H:=\max_{z\in\alpha}\delta(z)$, it follows from Lemma \ref{lowconvex} that
\[\log\left(\frac{H}{\sqrt{\delta(x)\delta(y)}}\right)\leq d_K(x,y)\leq\log \left(1+\frac{A|x-y|}{\sqrt{\delta(x) \delta(y)}}\right),\]
which implies $H\lesssim |x-y|+\sqrt{\delta(x)\delta(y)}$. As
$$
\begin{aligned}
\frac{L([x,y])}{H^{\frac{1}{m}}}&\leq \int_0^1\frac{|\dot{\alpha}(t)|}{\delta(\alpha(t))^{\frac{1}{m}}}dt\lesssim \int_0^1K(\alpha(t),\dot{\alpha}(t))dt=d_K(x,y)\\
&\leq\log\left(1+\frac{A|x-y|}{\sqrt{\delta_{\Omega}(x) \delta_{\Omega}(y)}}\right)\lesssim\frac{|x-y|}{\sqrt{\delta(x) \delta(y)}},
\end{aligned}
$$
we obtain
\[
L([x,y])\lesssim\frac{\left(|x-y|+\sqrt{\delta(x)\delta(y)}\right)^{\frac{1}{m}}|x-y|}{\sqrt{\delta(x) \delta(y)}}\lesssim|x-y|^{\frac{1}{m}}.
\]

If $(\Omega,\: d_K)$ is Gromov hyperbolic with $\delta>0$, then it follows from Theorem \ref{stab} that there exists a constant $R = R(\delta,\lambda)$ such that the Hausdorff distance between $[p, q]$ and the image of $\gamma$ is less than $R$. Denoting $H_\gamma:=\max_{\omega\in\gamma}\delta(\omega)=\delta(\omega_0)$, then there exists a point $z_0\in[x,y]$ with $d_K(\omega_0,z_0)\leq R$.

By Lemma \ref{lowconvex}, we have $\delta(\omega_0)\asymp\delta(z_0)$, which implies that \[H_\gamma\lesssim\delta(z_0)\leq H\lesssim|x-y|+\sqrt{\delta(x)\delta(y)}.\]
It follows from Corollary \ref{mco2} and $d_K(\omega_0,z_0)\leq R$ that
\[L([w_0,z_0])^m\lesssim \sqrt{\delta(w_0)\delta(z_0)}\lesssim H\lesssim|x-y|+\sqrt{\delta(x)\delta(y)}.\]
Therefore, if $\sqrt{\delta(x)\delta(y)}\leq|x-y|$, then we obtain $L([w_0,z_0])\lesssim |x-y|^{\frac{1}{m}}$ and
\[L(\gamma)\leq L([x,y])+2L([w_0,z_0])\lesssim|x-y|^{\frac{1}{m}}.\]
And when $\sqrt{\delta(x)\delta(y)}>|x-y|$, similar to the case of geodesic we have
\[
L(\gamma)\leq\lambda H_\gamma^{\frac{1}{m}}d_K(x,y)\lesssim\frac{\left(|x-y|+\sqrt{\delta(x)\delta(y)}\right)^{\frac{1}{m}}|x-y|}{\sqrt{\delta(x) \delta(y)}}\lesssim|x-y|^{\frac{1}{m}},
\]
which completes the proof.
$\hfill\qed$

\begin{thm}(=Theorem \ref{add})
Let $\Omega\subset \mathbb{C}^n(n\geq2)$ be a bounded domain. Suppose that $\Omega\cap U$ is $m$-convex and $\partial\Omega\cap U$ is Dini-smooth in a neighborhood $U$ of $\xi\in\partial\Omega$.
Then there exist a neighborhood $V$ of $\xi$ with $V\subset\subset U$ and a constant $C>0$ with
\begin{align}\label{loc11}
 d_{K_\Omega}(x,y)\leq d_{K_{\Omega\cap U}}(x,y)\leq d_{K_\Omega}(x,y)+C|x-y|^{\frac{1}{m}}\text{ for }x,y\in\Omega\cap V,
\end{align}
and for any $x\neq y\in\Omega\cap V$,
\begin{align}\label{loc22}
1\leq\frac{d_{K_{\Omega\cap U}}(x,y)}{d_{K_\Omega}(x,y)}\leq 1+C\left(\delta_\Omega(x)\wedge\delta_\Omega(y)+|x-y|^{\frac{1}{m}}\right).
\end{align}
\end{thm}

\noindent$Proof$. By \cite[Corollary 6.14]{22Visi}, it follows that there exist a neighborhood $V_0$ of $\xi$ with $V_0\subset\subset U$ and a constant $C>0$ such that, for $z\in\Omega\cap V_0$ and $X\in\mathbb{C}^n$,
\begin{align}\label{loc}
K_\Omega(z,X)\leq K_{\Omega\cap U}(z,X)\leq\left(1+C\delta_\Omega(z)\right)K_\Omega(z,X).
\end{align}

Then we only need to check the right side of inequality (\ref{loc11}) and (\ref{loc22}). By applying \cite[Proposition 4.4]{BZ22} and \cite[Lemma 16.5]{Zimmer2016Gromov}, we can choose a neighborhood $V$ of $\xi$ with $V\subset V_0$ such that, for any $x,y\in\Omega\cap V$ and $\varepsilon>0$, there exists a $(1,\varepsilon)$-quasi-geodesic $\gamma_\varepsilon$ of $d_{K_\Omega}$ connecting $x$ and $y$ which lies in $\Omega\cap V_0$.

Then by (\ref{loc}) and the lemmas in Section \ref{convex}, we deduce that $d_{K_\Omega}$ satisfies (\ref{slow}) and (\ref{fmet}) with $\alpha=\frac{1}{2},\beta=\frac{1}{m}$ in $\Omega\cap V$. By using the proof of Theorem \ref{Fin} and the local version of Lemma \ref{up} (see \cite[Theorem 7]{2015Estimates}), it follows that $L(\gamma_\varepsilon)\lesssim|x-y|^\frac{1}{m}+\varepsilon$.

For $z\in\Omega\cap V$ and $X\in\mathbb{C}^n$, as $K_\Omega(z,X)\leq\frac{|X|}{\delta_\Omega(z)}$, then
\[K_{\Omega\cap U}(z,X)\leq K_\Omega(z,X)+C|X|.\]
Then for $x,y\in\Omega\cap V$, by considering $(1,\varepsilon)$-quasi-geodesic $\gamma_\varepsilon$, we now obtain that
\[
d_{K_{\Omega\cap U}}(x,y)\leq L_{K}(\gamma_\varepsilon)+CL(\gamma_\varepsilon)\leq L_{K}(\gamma_\varepsilon)+C\left(|x-y|^{\frac{1}{m}}+\varepsilon\right).
\]
In particular, we obtain when $\varepsilon \rightarrow 0$
\[
d_{K_{\Omega\cap U}}(x,y)\leq d_{K_\Omega}(x,y)+C|x-y|^{\frac{1}{m}}.
\]
In addition, for $z\in\gamma_\varepsilon$ it follows that
$$
\delta_\Omega(z)\leq\delta_\Omega(x)\wedge\delta_\Omega(y)+L(\gamma_\varepsilon)
\lesssim\delta_\Omega(x)\wedge\delta_\Omega(y)+|x-y|^{\frac{1}{m}}+\varepsilon.
$$
Therefore, for $x\neq y\in\Omega\cap V$, we deduce that
\[
\frac{d_{K_{\Omega\cap U}}(x,y)}{d_{K_\Omega}(x,y)}\leq 1+C\left(\delta_\Omega(x)\wedge\delta_\Omega(y)+|x-y|^{\frac{1}{m}}\right),
\]
which completes the proof.
$\hfill\qed$

\medskip

Now we give the estimate for the Catlin-type distance $d_{\widetilde{K}}$ by using Theorem \ref{Fin} and the techniques in \cite[Proposition 2.5]{1987Localization}.

\begin{thm}
Suppose that $\Omega \subset \mathbb{C}^{2}$ is a smoothly bounded pseudoconvex domain of finite type $m$. Then there exists a constant $C \geq 1$ such that
$$
 g_m(x,y)-C\leq d_{\widetilde{K}}(x,y) \leq g_1(x,y)+C
$$
for each $x,y \in \Omega$.
\end{thm}

\noindent$Proof$.
First we estimate the upper bounds of $d_{\widetilde{K}}$. Note that we only need to consider the case of $x,y \in \Omega \cap U$ for a small neighborhood $U$ of a point $\xi \in\partial\Omega$. Denote by $p=\pi(x)$ and $q=\pi(y)$ the projections of $x$ and $y$ to the boundary and denote $x^{\prime}=x-|x-y|\vec{n}(p), \: y^{\prime}=y-|x-y|\vec{n}(q)$.

Then Lemma \ref{low} implies that
$$
d_{\widetilde{K}}\left(x, x^{\prime}\right)= \log \left(\frac{\delta\left(x^{\prime}\right)}{\delta(x)}\right)=\log \left(1+\frac{|x-y|}{\delta(x)}\right),
$$
and
$$
d_{\widetilde{K}}\left(y, y^{\prime}\right)= \log \left(\frac{\delta\left(y^{\prime}\right)}{\delta(y)}\right)=\log \left(1+\frac{|x-y|}{\delta(y)}\right).
$$

It remains to find an upper bound for $d_{\widetilde{K}}\left(x^{\prime}, y^{\prime}\right)$.
By Lemma \ref{pro}, for a sufficiently small neighborhood $U$, we obtain that $|\vec{n}(p)-\vec{n}(q)|\leq\varepsilon<1$ for $p,\: q\in\partial\Omega\cap U$. It follows that
$$|x^{\prime}-y^{\prime}|\leq(1+|\vec{n}(p)-\vec{n}(q)|)|x-y|\leq(1+\varepsilon)|x-y|<2|x-y|,$$
and
$$\delta(x^{\prime})\wedge\delta(y^{\prime})\geq\delta(x)\wedge\delta(y)+|x-y|>|x-y|.$$

Consider the mapping $\Phi:\mathbb{C}\rightarrow\mathbb{C}^n$ defined by
$\Phi(\xi)=x^{\prime}+\xi(y^{\prime}-x^{\prime})$.
Then, writing
$$D:=\left\{\xi\in\mathbb{C}:\min\left\{|\xi|,|\xi-1|\right\}<\frac{1}{1+\varepsilon}\right\},$$
we have $\Phi(D)\subset\Omega\cap U$ with $x^{\prime}=\Phi(0),\: y^{\prime}=\Phi(1)$.
Then
$$d_{\widetilde{K}}\left(x^{\prime}, y^{\prime}\right)\lesssim d_{K_\Omega}\left(x^{\prime}, y^{\prime}\right)
\leq d_{K_D}\left(0, 1\right)\leq C.$$
As a consequence,
$$
\begin{aligned}
d_{\widetilde{K}}(x, y) &\leq d_{\widetilde{K}}\left(x, x^{\prime}\right)+d_{\widetilde{K}}\left(x^{\prime}, y^{\prime}\right)+d_{\widetilde{K}}\left(y, y^{\prime}\right)\\
&\leq \log \left(1+\frac{|x-y|}{\delta(x)}\right)+\log \left(1+\frac{|x-y|}{\delta(y)}\right)+C\\
&\leq \log \left(\frac{(\delta(x)+|x-y|)(\delta(y)+|x-y|)}{\delta(x)\delta(y)}\right)+C\\
&\leq 2\log \left(\frac{|x-y|+\delta(x)\vee\delta(y)}{\sqrt{\delta(x)\delta(y)}}\right)+C\\
&\leq g_1(x,y)+C,
\end{aligned}
$$
which proves the proof in this case.

Now we estimate the lower bounds of $d_{\widetilde{K}}$ by Theorem \ref{Fin}.
Choose a sufficiently small $\varepsilon$ such that we can use Lemma \ref{low} and the inequality (\ref{nk}) on $N_\varepsilon(\partial\Omega)$. We only need to prove the case of $x,y \in \Omega \cap N_\varepsilon(\partial\Omega)$.

Lemma \ref{low} implies that $d_{\widetilde{K}}$ satisfied (\ref{slow}) with $\alpha=1$. Moreover, by the inequality (\ref{nk}), for $z\in\Omega\cap N_\varepsilon(\partial\Omega)$ and $X\in\mathbb{C}^2$, it follows that
\[
\widetilde{K}(z,X)\geq\frac{|X_N|}{\delta(z)}+\frac{1}{C}\frac{|X|}{\delta(z)^{\frac{1}{m}}}
\gtrsim\frac{|X|}{\delta(z)^{\frac{1}{m}}},
\]
which satisfies (\ref{fmet}).
Suppose $\gamma:[0,1] \rightarrow \Omega$ is a piecewise $C^1$-smooth curve with endpoints $x=\gamma(0)$ and $y=\gamma(1)$. If $\gamma\subset N_\varepsilon(\partial\Omega)$, then by the argument in the proof of Theorem \ref{Fin}, we obtain that
$$
\begin{aligned}
\widetilde{K} \text {-length }(\gamma)&\geq 2\log \left(\frac{L(\gamma)^m+\delta(x)\vee\delta(y)}{\sqrt{\delta(x)\delta(y)}}\right)-C\\
&\geq 2\log \left(\frac{|x-y|^m+\delta(x)\vee\delta(y)}{\sqrt{\delta(x)\delta(y)}}\right)-C.
\end{aligned}
$$
If $\gamma$ is not in $N_\varepsilon(\partial\Omega)$, the we set
\[t_0:=\min\{t\in[0,1]:\gamma(t)\notin N_\varepsilon(\partial\Omega)\},\]
and
\[t_1:=\max\{t\in[0,1]:\gamma(t)\notin N_\varepsilon(\partial\Omega)\}.\]
Then $\gamma_0:=\gamma|_{[0,t_0]}$ and $\gamma_1:=\gamma|_{[t_1,1]}$ in $N_\varepsilon(\partial\Omega)$. Lemma \ref{low} implies that
\[\widetilde{K}\text {-length }(\gamma_0)\geq \log \left(\frac{\varepsilon}{\delta(x)}\right)\text{ and }\widetilde{K}\text {-length }(\gamma_0)\geq \log \left(\frac{\varepsilon}{\delta(y)}\right).\]
Since there exists a constant $M>0$ with $|x-y|^m\leq M\varepsilon$ for all $x,y\in\Omega$, we now obtain
\[\widetilde{K}\text {-length }(\gamma)\geq 2\log \left(\frac{\varepsilon}{\sqrt{\delta(x)\delta(y)}}\right)\geq 2\log \left(\frac{|x-y|^m}{\sqrt{\delta(x)\delta(y)}}\right)-C.\]
As $\gamma$ is arbitrary, we obtain from Lemma \ref{low} that
$$
\begin{aligned}
d_{\widetilde{K}}(x, y) &\geq 2\log \left(\frac{|x-y|^m+\delta(x)\vee\delta(y)}{\sqrt{\delta(x)\delta(y)}}\right)-C\\
&\geq g_m(x,y)-C,
\end{aligned}
$$
which completes the proof.
$\hfill\qed$

\bigskip
\section{Boundary Extensions of Mappings}\label{exten}

Now we are ready to prove Theorem \ref{miso}. For the convenience of reader, we restate it as follows:
\begin{thm}
 Let $\Omega_1,\: \Omega_2 \subset \mathbb{C}^n(n\geq2)$ be bounded $m$-convex domains with Dini-smooth boundary, and let $\overline{\Omega}_1,\overline{\Omega}_2$ be their Euclidean compactifications respectively. If $f: \Omega_1 \rightarrow \Omega_2$ is an isometry with respect to the Kobayashi metrics, then $f$ has a homeomorphic extension $\bar{f}: \overline{\Omega}_1 \rightarrow \overline{\Omega}_2$ such that the induced boundary map $\bar{f}|_{\partial \Omega_1}: \partial \Omega_1 \rightarrow \partial \Omega_2$ is bi-H\"{o}lder continuous with respect to the Euclidean metric.
\end{thm}

\noindent$Proof$. From Lemma \ref{lowconvex} and Lemma \ref{up}, it follows that there exists a constant $M$ such that, for all $x,y\in\Omega_1$,
\[
\frac{1}{2}\log \left(\frac{\delta_{\Omega_1}(y)}{\delta_{\Omega_1}(x)}\right)\leq d_{K_{\Omega_1}}(x,y)\leq\frac{1}{2}\log \left(\frac{M}{\delta_{\Omega_1}(x)\delta_{\Omega_1}(y)}\right),
\]
and
\[
\frac{1}{2}\log \left(\frac{\delta_{\Omega_2}(f(y))}{\delta_{\Omega_2}(f(x))}\right)\leq d_{K_{\Omega_2}}(f(x),f(y))\leq\frac{1}{2}\log \left(\frac{M}{\delta_{\Omega_2}(f(x))\delta_{\Omega_2}(f(y))}\right).
\]
 Noting that $f$ is an isometry, then
\[
\frac{1}{M}\delta_{\Omega_2}(f(y))\delta_{\Omega_1}(y)\delta_{\Omega_1}(x)\leq\delta_{\Omega_2}(f(x))
\leq\frac{M\delta_{\Omega_2}(f(y))}{\delta_{\Omega_1}(y)}\delta_{\Omega_1}(x),
\]
which implies $\delta_{\Omega_2}(f(x))\asymp\delta_{\Omega_1}(x)$ by fixing $y=y_0$.

Now it follows from Corollary \ref{mco} that there exists a constant $C>0$ such that, for all $x,y\in\Omega_1$,
\[
\left|x-y\right|^{m}\lesssim\left|f(x)-f(y)\right|+\sqrt{\delta_{\Omega_2}(f(x))\delta_{\Omega_2}(f(y))},
\]
and
\[
\left|f(x)-f(y)\right|^{m}\lesssim\left|x-y\right|+\sqrt{\delta_{\Omega_1}(x)\delta_{\Omega_1}(y)}.
\]

Then take two sequences $\{x_k\}$ and $\{y_k\}$ in $\Omega_1$ with $x_k\rightarrow a\in\partial\Omega_1$ and $y_k\rightarrow b\in\partial\Omega_1$.
Noting that $\sqrt{\delta_{\Omega_1}(x_k)\delta_{\Omega_1}(y_k)}\rightarrow0$ and $\sqrt{\delta_{\Omega_2}(f(x_k))\delta_{\Omega_2}(f(y_k))}\rightarrow0$,
it is easy to see that $\left|f(x_k)-f(y_k)\right|\rightarrow0$ if $a=b$, which means $f$ extends continuously to the boundary.

Moreover, for $a,b\in \partial\Omega_1$, we therefore obtain the inequalities
\[
\left|a-b\right|^{m}\lesssim\left|f(a)-f(b)\right|\text{ and }\left|f(a)-f(b)\right|^{m}\lesssim\left|a-b\right|,
\]
which imply the following bi-H\"{o}lder estimates
\[|a-b|^{m}\lesssim|f(a)-f(b)|\lesssim|a-b|^{\frac{1}{m}}.\]
It means that $f$ not only extends to a homeomorphism on $\overline{\Omega}_1$ but also bi-H\"{o}lder continuous on $\partial\Omega_1$.
The proof is now complete.
$\hfill\qed$

\bigskip
We can now prove that the identity map for the smooth bounded pseudoconvex domain of finite type in $\mathbb{C}^2$ can extend to a bi-H\"{o}lder map between the Euclidean boundary and Gromov boundary.
\begin{thm}
Suppose that $\Omega \subset \mathbb{C}^{2}$ is a bounded pseudoconvex domain with smooth boundary of finite D'Angelo type $m$. Then the identity map id : $\Omega \rightarrow \Omega$ can extend to a bi-H\"{o}lder map
$$
id:(\partial \Omega ,|\cdot|) \rightarrow\left(\partial_G^{d_{K}}\Omega, \rho_G^{d_{K}}\right)
$$
between the boundaries, where $\rho_G^{d_{K}}$ is a visual metric on the Gromov boundary of $\left(\Omega, d_{K}\right)$.
\end{thm}

$Proof.$ Fix a base point $\omega\in\Omega$. For any $x,\: y\in\Omega$, it follows from Theorem \ref{dg} that there exist constants $C_1,C_2>0$ with
\begin{align}\label{inequ}
\log\frac{C_1}{\left|x-y\right|+\delta(x)\vee\delta(y)}\leq(x|y)_\omega^{d_{\widetilde{K}}}
\leq\log\frac{C_2}{\left|x-y\right|^{m}+\delta(x)\vee\delta(y)}.
\end{align}

Then for any Gromov sequences $\{x_k\},\{y_k\}\subset\Omega$ with $(x_k|y_k)_\omega^{d_{\widetilde{K}}}\rightarrow\infty$ as $k\rightarrow\infty$,
using the right side of (\ref{inequ}), we get $x_k,y_k$ tend to the boundary and $\left|x_k-y_k\right|\rightarrow0$ as $k\rightarrow\infty$.
On the other hand, for every point $a\in\partial\Omega$, choose sequences $\{x_k\},\{y_k\}\subset\Omega$ with $x_k\rightarrow a$ and $y_k\rightarrow a$ as $k\rightarrow\infty$. By the left side of (\ref{inequ}), we get $(x_k|y_k)_\omega^{d_{\widetilde{K}}}\rightarrow\infty$ as $k\rightarrow\infty$. It means a sequence in $\Omega$ is a Gromov sequence with $d_{\widetilde{K}}$ if and only if it converges to certain boundary point in the Euclidean metric. So the identity map can extend to a bijection between the Euclidean boundary of $\Omega$ and the Gromov boundary of the space $\left(\Omega, d_{\widetilde{K}}\right)$.

Moreover, for $a,b\in \partial\Omega$, letting $x$ and $y$ tend to $a$ and $b$ respectively in (\ref{inequ}), we now obtain that
\[
|a-b|^{m}\lesssim\exp\left(-\left(a|b\right)_\omega^{d_{\widetilde{K}}}\right)\lesssim|a-b|.
\]
Note that, for visual metrics $\rho_G^{d_{\widetilde{K}}}$ on $\partial_G^{d_{\widetilde{K}}}\Omega$, there exists a parameter $\epsilon>0$ such that
$$
\rho_G^{d_{\widetilde{K}}}(a, b) \asymp \exp \left(-\epsilon(a|b)_\omega^{d_{\widetilde{K}}}\right), \text { for } a, b \in \partial_G^{d_{\widetilde{K}}} \Omega.
$$
It means the boundary mapping $id:(\partial \Omega ,|\cdot|) \rightarrow\left(\partial_G^{d_{\widetilde{K}}}\Omega, \rho_G^{d_{\widetilde{K}}}\right)$ is bi-H\"{o}lder continuous. Thus, we complete the proof by using Proposition \ref{visual}.
$\hfill\qed$

\bigskip
Finally we conclude this section by proving Theorem \ref{bi}. For the convenience of reader, we restate it as follows:
\begin{thm}
 Let $\Omega_1,\Omega_2 \subset \mathbb{C}^2$ be bounded smooth pseudoconvex domains of finite D'Angelo type with smooth boundary, and let $\overline{\Omega}_1,\overline{\Omega}_2$ be their Euclidean compactifications. If $f: \Omega_1 \rightarrow \Omega_2$ is a quasi-isometry with respect to the Kobayashi metrics, then $f$ has a homeomorphic extension $\bar{f}: \overline{\Omega}_1 \rightarrow \overline{\Omega}_2$ such that the induced boundary map $\bar{f}|_{\partial \Omega_1}: \partial \Omega_1 \rightarrow \partial \Omega_2$ is bi-H\"{o}lder continuous with respect to the Euclidean metric.
\end{thm}

$Proof.$ Note the Gromov hyperbolicity of $\Omega_1,\Omega_2$ endowed with Kobayashi metrics. In view of Proposition \ref{extend}, the quasi-isometry $f$ induces a bi-H\"{o}lder mapping
$\tilde{f}: \partial_G^{d_{K_1}}\Omega_1\rightarrow \partial_G^{d_{K_2}}\Omega_2,$
where the Gromov boundaries $\partial_G^{d_{K_i}}\Omega_i$ of $\left(\Omega_i,d _{K_i}\right), \: i=1,2$, are endowed with their visual metrics.

Recall that $\rho_G^{d_{K_i}}$ is the visual metric of the Gromov boundary of $\left(\Omega_i, d _{K_i}\right)$. It follows from Theorem \ref{id} that the identity map $i d_i: \Omega_i \rightarrow \Omega_i, \:i=1,2$, induces a bi-H\"{o}lder map
$$
i d_i:\left(\partial \Omega_i,|\cdot|\right) \rightarrow\left(\partial_G^{d_{K_i}}\Omega_i, \rho_G^{d_{K_i}}\right).
$$
Therefore,
$$
\bar{f}=i d_2^{-1} \circ \tilde{f} \circ i d_1: \partial \Omega_1 \rightarrow \partial\Omega_2.
$$
is a well-defined boundary map.

Noting that $\bar{f}$ is the composition of bi-H\"{o}lder maps, it is also a bi-H\"{o}lder map. The proof is now complete.
$\hfill\qed$

\bigskip
{\bf Funding}. Jinsong Liu is supported by National Key R\&D Program of China (Grant No. 2021YFA1003100), NSFC (Grants No. 11925107, and No. 12226334), Key Research Program of Frontier Sciences, CAS (Grant No. ZDBS-LY-7002). Hongyu Wang is supported by the Fundamental Research Funds for the Central Universities (500422379).

\vspace{0.3cm} \noindent{\bf Acknowledgements}. The authors
would like to thank the referee for a careful reading and valuable comments.

\bibliography{reference}
\bibliographystyle{plain}{}
\end{document}